\documentclass[12pt]{article}

\usepackage{amsmath,amsthm}
\usepackage{amssymb}
\usepackage{pstricks,pst-node}
\usepackage{color}
\usepackage{xspace}
\usepackage[colorlinks=true,
linkcolor=green,
filecolor=brown,
citecolor=green]{hyperref}

\def\v{\vert}
\def\a{\ensuremath{\mathcal A}\xspace}
\def\b{\ensuremath{\mathcal B}\xspace}
\def\c{\ensuremath{\mathcal C}\xspace}
\def\C{\ensuremath{\mathsf{C}}\xspace}
\def\si{\sigma}
\def\ep{\epsilon}
\def\el{\ensuremath{\mathcal L}\xspace}

\def\ins{nonfirst\xspace}

\newcommand{\StirlingCycle}[2]{\genfrac{[}{]}{0pt}{}{#1}{#2}}

\begin{document}
\newtheorem*{main}{Main Theorem}
\newtheorem{theorem}{Theorem}
\newtheorem{defn}[theorem]{Definition}
\newtheorem{lemma}[theorem]{Lemma}
\newtheorem{prop}[theorem]{Proposition}
\newtheorem{cor}[theorem]{Corollary}
\mbox{} 
\vspace*{-15mm}
\begin{center}
{\Large
A Determinant of Stirling Cycle Numbers 
Counts Unlabeled Acyclic Single-Source Automata                            \\ 
}

\vspace*{5mm}

DAVID CALLAN  \\
Department of Statistics  \\
\vspace*{-2mm}
University of Wisconsin-Madison  \\
\vspace*{-2mm}
1300 University Ave  \\
\vspace*{-2mm}
Madison, WI \ 53706-1532  \\
{\bf callan@stat.wisc.edu}  \\
\vspace{5mm}

March 30, 2007
\end{center}

\begin{abstract}
We show that a determinant of Stirling 
cycle numbers counts unlabeled acyclic single-source automata. The  
proof involves a bijection from these automata to certain marked 
lattice paths and a sign-reversing involution to evaluate the determinant.

\end{abstract}

\vspace{3mm}

{\Large \textbf{1 \quad  Introduction}  } The chief purpose of this 
paper is to show bijectively that a determinant of Stirling cycle 
numbers counts unlabeled acyclic single-source automata. 
Specifically, let 
 $A_{k}(n)$ denote the $kn \times kn$ matrix with 
$(i,j)$ entry $\StirlingCycle{\lfloor \frac{i-1}{k} \rfloor+2}{\lfloor 
\frac{i-1}{k}\rfloor +1+i-j}$, where $\StirlingCycle{i}{j}$ is the 
Stirling cycle 
number, the number of permutations on $[i]$ with $j$ cycles. For example,
\[ 
A_{2}(5)=
\left(\begin{array}{cccccccccc}
    1 & 0 & 0 & 0 & 0 & 0 & 0 & 0 & 0 & 0 \\ 
1 & 1 & 0 & 0 & 0 & 0 & 0 & 0 & 0 & 0 \\ 
0 & 1 & 3 & 2 & 0 & 0 & 0 & 0 & 0 & 0 \\
0 & 0 & 1 & 3 & 2 & 0 & 0 & 0 & 0 & 0 \\ 
0 & 0 & 0 & 1 & 6 & 11 & 6 & 0 & 0 & 0 \\ 
0 & 0 & 0 & 0 & 1 & 6 & 11 & 6 & 0 & 0 \\
0 & 0 & 0 & 0 & 0 & 1 & 10 & 35 & 50 & 24 \\ 
0 & 0 & 0 & 0 & 0 & 0 & 1 & 10 & 35 & 50 \\ 
0 & 0 & 0 & 0 & 0 & 0 & 0 & 1 & 15 & 85 \\
0 & 0 & 0 & 0 & 0 & 0 & 0 & 0 & 1 & 15 
\end{array}\right).
\]
As evident in the example, $A_{k}(n)$ is formed from $k$ copies of each of 
rows 2 through $n+1$ of the Stirling cycle triangle, arranged so that the first 
nonzero entry in each row is a 1 and, after the first row, this 1 
occurs just before the main diagonal; in other words, $A_{k}(n)$ is a Hessenberg 
matrix with 1s on the infra-diagonal. We will show 
\begin{main}
The determinant of $A_{k}(n)$ is the number of unlabeled acyclic single-source automata 
with $n$ transient states on a $(k+1)$-letter input alphabet. 
\end{main}
Section 2 reviews basic terminology 
for automata and recurrence relations to count finite acyclic automata. 
Section 3 introduces 
column-marked subdiagonal paths, which play an intermediate role, and a 
way to code them. Section 4
presents a bijection from these column-marked subdiagonal paths to 
unlabeled acyclic single-source automata. Finally, Section 5 evaluates 
$\det A_{k}(n)$ using a sign-reversing involution and shows that the 
determinant counts the codes for column-marked subdiagonal paths.
\vspace{10mm}

{\Large \textbf{2 \quad  Automata}  }

A (complete, deterministic) automaton consists of a set of states and 
an input alphabet whose letters transform the states among 
themselves: a letter and a state produce another state (possibly the 
same one). A finite automaton (finite set of states, finite input 
alphabet of, say, $k$ letters) can be represented as a $k$-regular 
directed multigraph with ordered edges: the vertices represent the 
states and the first, second, \ldots edge from a vertex give the 
effect of the first, second, \ldots alphabet letter on that state.
A finite automaton cannot be acyclic in the usual sense of no cycles: 
pick a vertex and follow any path from it. This path must ultimately 
hit a previously encountered vertex, thereby creating a cycle. So the 
term acyclic is used in the looser sense that only one vertex, called 
the \emph{sink}, is involved in cycles. This means that all edges from 
the sink loop back to itself (and may safely be omitted) 
and all other paths feed into the sink. 

A non-sink state is called \emph{transient}. The \emph{size} of an 
acyclic automaton is the number of transient states. An acyclic 
automaton of size $n$ thus has transient states which we label 
$1,2,\ldots,n$ and a sink, labeled $n+1$. Liskovets 
\cite{Liskovets03} uses the inclusion-exclusion principle (more about 
this below) to obtain the following recurrence relation for the 
number $a_{k}(n)$ of acyclic automata of size $n$ on a $k$-letter 
input alphabet ($k\ge 1$):
\[
a_{k}(0)=1;\qquad 
a_{k}(n)=\sum_{j=0}^{n-1}(-1)^{n-j-1}\binom{n}{j}(j+1)^{k(n-j)}a_{k}(j),\quad 
n\ge 1.
\]

A \emph{source} is a vertex with no incoming edges. A finite acyclic 
automaton has at least one source because a path traversed backward 
$v_{1} \leftarrow v_{2} \leftarrow v_{3} \leftarrow \ldots$ must have 
distinct vertices and so cannot continue indefinitely. An automaton 
is \emph{single-source} (or initially connected) if it has only one 
source. Let $\b_{k}(n)$ denote the set of single-source acyclic 
finite (SAF) automata on a $k$-letter input alphabet with vertices 
$1,2,\ldots,n+1$ where $1$ is the source and $n+1$ is the sink, and 
set $b_{k}(n)=\v\, \b_{k}(n) \,\v$. The \emph{two-line representation} of 
an automaton in $\b_{k}(n)$ is the $2 \times kn$ matrix whose columns 
list the edges in order. For example,
\[
B=
\left(\begin{array}{ccccccccccccccc}
    1 & 1 & 1 & 2 & 2 & 2 & 3 & 3 & 3 & 4 & 4 & 4 & 5 & 5 & 5  \\
    2 & 4 & 6 & 6 & 6 & 6 & 6 & 6 & 6 & 3 & 5 & 3 & 2 & 2 & 6
\end{array}
\right) 
\]
is in $\b_{3}(5)$ and the source-to-sink paths in $B$ include 
$1\overset{a}{\:\rightarrow} 2 \overset{a}{\:\rightarrow} 6,\ 
1\overset{b}{\:\rightarrow} 4 \overset{c}{\:\rightarrow} 
3 \overset{a}{\:\rightarrow} 
6,\ 1 \overset{b}{\:\rightarrow} 4 \overset{b}{\:\rightarrow} 
5 \overset{b}{\:\rightarrow} 2 \overset{b}{\:\rightarrow} 6$, where 
the alphabet is $\{a,b,c\}$.
\begin{prop}\label{recurrence}
    The number $b_{k}(n)$ of SAF automata of size $n$ on a $k$-letter 
    input alphabet $(n,k\ge 1)$ is given by 
    \[
     b_{k}(n)=\sum_{i=1}^{n}(-1)^{n-i}
    \binom{n-1}{i-1}(i+1)^{k(n-i)}a_{k}(i)
    \]
\end{prop}

\textbf{Remark}\quad This formula is a bit more succinct than the the 
recurrence in \cite[Theorem 3.2]{Liskovets03}.

\textbf{Proof}\quad Consider the set $\a$ of acyclic automata with 
transient vertices $[n]=\{1,2,\ldots,n\}$ in which 1 is a source. 
Call $2,3,\ldots,n$ the \emph{interior} vertices. For $X\subseteq 
[2,n]$, let
\begin{eqnarray*}
   f(X)   & = &  \textrm{\# automata in \a whose set of interior vertices 
includes $X$}, \\
    g(X) & = & \textrm{\# automata in \a whose set of interior vertices 
is precisely $X$}.
\end{eqnarray*}
Then $f(X)=\sum_{Y:\,X\subseteq Y \subseteq [2,n]}g(Y)$ and by 
M\"{o}bius inversion \cite{lintwilson} on the lattice of subsets of 
$[2,n],\ g(X)=\sum_{Y:\,X\subseteq Y \subseteq [2,n]}\mu(X,Y)f(Y)$ 
where $\mu(X,Y)$ is the M\"{o}bius function for this lattice. Since 
$\mu(X,Y)=(-1)^{\v Y \v-\v X \v}$ if $X\subseteq Y$, we have in 
particular that 
\begin{equation}
g(\emptyset)=\sum_{Y \subseteq [2,n]}(-1)^{\v\, Y \,\v}f(Y).    
    \label{eq:dagger}
\end{equation}
Let $\v\, Y \,\v=n-i$ so that $1\le i \le n$. When $Y$  consists 
entirely of sources, the vertices in $[n+1]\backslash Y$ and their 
incident edges form a subautomaton with $i$ transient states; there 
are $a_{k}(i)$ such. Also, all edges from the $n-i$ 
vertices comprising $Y$ go directly into $[n+1]\backslash Y:\ 
(i+1)^{k(n-i)}$ choices. Thus $f(Y)=(i+1)^{k(n-i)}a_{k}(i)$. By 
definition, $g(\emptyset)$ is the number of automata in \a for which 
1 is the only source, that is,  $g(\emptyset) = b_{k}(n)$ and the 
Proposition now follows from (\ref{eq:dagger}). \qed

An \emph{unlabeled} SAF automaton is an equivalence class of SAF automata 
under relabeling of the interior vertices. Liskovets notes 
\cite{Liskovets03} (and we prove below) that $\b_{k}(n)$ has no 
nontrivial automorphisms, that is, each of the $(n-1)!$ relabelings 
of the interior vertices of $B \in \b_{k}(n)$ produces a different 
automaton. So unlabeled SAF automata of size $n$ on a $k$-letter 
alphabet are counted by 
$\frac{1}{(n-1)!}b_{k}(n)$. The next result establishes a 
canonical representative in each relabeling class.
\begin{prop}
    Each equivalence class in $\b_{k}(n)$ under relabeling of 
    interior vertices has size $(n-1)!$ and contains exactly one SAF 
    automaton with the ``last occurrences increasing'' property: the last 
    occurrences of the interior vertices---$2,3,\ldots,n$---in the 
    bottom row of its two-line representation occur in that order.
\end{prop}
\textbf{Proof}\quad The first assertion follows from the fact that the 
interior vertices of an automaton $B\in b_{k}(n)$ can be distinguished 
intrinsically, that is, independent of their labeling. To see this, 
first mark the source, namely 1, with a mark (new label) $v_{1}$ and 
observe that there exists at least one interior vertex whose only 
incoming edge(s) are from the source (the only currently marked vertex) 
for otherwise a cycle would be present. For each such interior vertex 
$v$, choose the last edge from the marked vertex to $v$ using the 
built-in ordering of these edges. This determines an order on these 
vertices; mark them in order $v_{2},v_{3},\ldots,v_{j}\ (j\ge 2)$. If there 
still remain unmarked interior vertices, at least one of them has incoming 
edges only from a marked vertex or again a cycle would be present. For 
each such vertex, use the last incoming edge from a marked vertex, 
where now edges are arranged in order of initial vertex $v_{i}$ with 
the built-in order breaking ties, to order and mark these vertices 
$v_{j+1},v_{j+2},\ldots $. Proceed similarly until all interior vertices 
are marked. For example, for 
\[
B=
\left(\begin{array}{ccccccccccccccc}
    1 & 1 & 1 & 2 & 2 & 2 & 3 & 3 & 3 & 4 & 4 & 4 & 5 & 5 & 5  \\
    2 & 4 & 6 & 6 & 6 & 6 & 6 & 6 & 6 & 3 & 5 & 3 & 2 & 2 & 6
\end{array}
\right),
\]
$v_{1}=1$ and there is just one interior vertex, namely 4, whose only 
incoming edge is from the source, and so $v_{2}=4$ and 4 becomes a 
marked vertex. Now all incoming edges to both 3 and 5 are from marked 
vertices and the last such edges (built-in order comes into play) are 
$4 \overset{b}{\:\rightarrow} 
5$ and $4 \overset{c}{\:\rightarrow} 3$ putting vertices 3,\,5 in the 
order 5,\,3. So $v_{3}=5$ and $v_{4}=3$. Finally, $v_{5}=2$.
This proves the first assertion. By construction of the $v$s, 
relabeling each interior vertex $i$ with the subscript of its 
corresponding $v$ produces an automaton in $\b_{k}(n)$ with the 
``last occurrences increasing'' property and is the only relabeling 
that does so. The example yields
\[
\left(\begin{array}{ccccccccccccccc}
    1 & 1 & 1 & 2 & 2 & 2 & 3 & 3 & 3 & 4 & 4 & 4 & 5 & 5 & 5  \\
    5 & 2 & 6 & 4 & 3 & 4 & 5 & 5 & 6 & 6 & 6 & 6 & 6 & 6 & 6
\end{array}
\right) .
\]
\qed

Now let $\c_{k}(n)$ denote the set of canonical SAF automata in 
$\b_{k}(n)$ representing unlabeled automata; thus $\v\,\c_{k}(n)\,\v=\frac{1}{(n-1)!}b_{k}(n)$.
Henceforth, we identify an unlabeled automaton with its canonical 
representative.

\vspace*{10mm}

{\Large \textbf{3 \quad  Column-Marked Subdiagonal Paths}  }

A \emph{subdiagonal} $(k,n,p)$-path is a lattice path of steps 
$E=(1,0)$ and $N=(0,1)$, $E$ for east and $N$ for north, from $(0,0)$ 
to $(kn,p)$ that never rise above the line $y=\frac{1}{k}x$. Let 
$\C_{k}(n,p)$ denote the set of such paths.For $k\ge 1$, it is clear 
that $\C_{k}(n,p)$ is nonempty only for $0\le p \le n$ and it is 
known (generalized ballot theorem) that 
\[
\v\, \C_{k}(n,p) \,\v = \frac{kn-kp+1}{kn+p+1}\binom{kn+p+1}{p}.
\]

A path $P$ in $\C_{k}(n,n)$ can be coded by the heights of its $E$ steps 
above the line $y=-1$; this gives a a sequence $(b_{i})_{i=1}^{kn}$ 
subject to the restrictions $1\le b_{1}\le b_{2}\le \ldots \le b_{kn}$ 
and $b_{i} \le \lceil i/k \rceil$ for all $i$.

A \emph{column-marked} subdiagonal $(k,n,p)$-path is one in which, 
for each $i \in [1,kn]$, one of the lattice squares below the $i$th 
$E$ step and above the horizontal line $y=-1$ is marked, say with a 
`$\,*\,$'.
Let $\C_{k}^{\textrm{\raisebox{1ex}{$*$}}}(n,p)$ denote the set 
of such  marked paths.
\begin{center}

\begin{pspicture}(-4,-1.5)(4,5.5)

\psline(-4,1)(-2,1)(-2,2)(1,2)(1,3)(4,3)(4,4)
\psline(-4,1)(4,5)
\psline[linecolor=red,linestyle=dotted](-4,0)(4,0)
\psline[linecolor=red,linestyle=dotted](-2,1)(4,1)
\psline[linecolor=red,linestyle=dotted](1,2)(4,2)
\psline[linecolor=red,linestyle=dotted](-4,0)(-4,1)
\psline[linecolor=red,linestyle=dotted](-3,0)(-3,1)
\psline[linecolor=red,linestyle=dotted](-2,0)(-2,1)
\psline[linecolor=red,linestyle=dotted](-1,0)(-1,2)
\psline[linecolor=red,linestyle=dotted](0,0)(0,2)
\psline[linecolor=red,linestyle=dotted](1,0)(1,2)
\psline[linecolor=red,linestyle=dotted](2,0)(2,3)
\psline[linecolor=red,linestyle=dotted](3,0)(3,3)
\psline[linecolor=red,linestyle=dotted](4,0)(4,3)

\psdots(-4,1)(-3,1)(-2,1)(-2,2)(-1,2)(0,2)(1,2)(1,3)(2,3)(3,3)(4,3)(4,4)
\psdots(0,3)(2,4)(4,5)

\rput(-3.5,0.5){\textrm{{\footnotesize $*$}}}
\rput(-2.5,0.5){\textrm{{\footnotesize $*$}}}
\rput(-1.5,0.5){\textrm{{\footnotesize $*$}}}
\rput(-0.5,1.5){\textrm{{\footnotesize $*$}}}
\rput(0.5,0.5){\textrm{{\footnotesize $*$}}}
\rput(1.5,2.5){\textrm{{\footnotesize $*$}}}
\rput(2.5,0.5){\textrm{{\footnotesize $*$}}}
\rput(3.5,1.5){\textrm{{\footnotesize $*$}}}

\rput(-4.5,1.2){\textrm{{\footnotesize (0,0)}}}
\rput(4.5,5.2){\textrm{{\footnotesize (8,4)}}}
\rput(0,-0.3){\textrm{{\footnotesize $y=-1$}}}
\rput(-1,3.0){\textrm{{\footnotesize $y=\frac{1}{2}x$}}}

\rput(0,-1.3){\textrm{{\small A path in $\C_{2}^{\textrm{\raisebox{1ex}{$*$}}}(4,3)$}}}

\end{pspicture}
\end{center} 

A marked path $P^{*}$ in $\C_{k}^{\textrm{\raisebox{1ex}{$*$}}}(n,n)$ can be coded by
a sequence of pairs $\big( (a_{i},b_{i})\big)_{i=1}^{kn}$ where 
$(b_{i})_{i=1}^{kn}$ is the code for the underlying path $P$ and 
$a_{i} \in [1,b_{i}]$ gives the position of the $*$ in the $i$th column. 
The example is coded by 
$(1,1),\,(1,1),\,(1,2),\,(2,2),\,(1,2),(3,3),\,(1,3),\,(2,3)$.

An explicit sum for $\v\, \C_{k}^{\textrm{\raisebox{1ex}{$*$}}}(n,n) 
\,\v $
is
\[
\v\, \C_{k}^{\textrm{\raisebox{1ex}{$*$}}}(n,n) \,\v =
\sum_{\substack{1\le b_{1}\le b_{2}\le \ldots \le b_{kn}, \\
\textrm{\raisebox{-1.5ex}{ $b_{i} \le \lceil i/k \rceil$  for all  $i$ } }}} b_{1} b_{2} 
\ldots b_{kn},
\]
because the summand $b_{1} b_{2} \ldots b_{kn}$ is the number of ways 
to insert the `$\,*\,$'s in the underlying path coded by 
$(b_{i})_{i=1}^{kn}$. 

It is also possible to obtain a recurrence for 
$\v\, \C_{k}^{\textrm{\raisebox{1ex}{$*$}}}(n,p) \,\v$, and then, using Prop. \ref{recurrence}, to 
show analytically that $\v\, \C_{k}^{\textrm{\raisebox{1ex}{$*$}}}(n,n) \,\v=
\v\, \c_{k+1}(n) \,\v$. 
However, it is much more pleasant to give a bijection 
and in the next section we will do so.
In particular, the number of SAF automata on 
a 2-letter alphabet is 
\[
\v\,\c_{2}(n)\,\v = \v\, 
\C_{1}^{\textrm{\raisebox{1ex}{$*$}}}(n,n) \,\v =
\sum_{\substack{1\le b_{1}\le b_{2}\le \ldots \le b_{n} \\
\textrm{\raisebox{-1.5ex}{ $b_{i} \le i$  for all  $i$ } }}} b_{1} b_{2} 
\ldots b_{n} = (1,3,16,127,1363,\ldots)_{n\ge 1},
\]
sequence
\htmladdnormallink{A082161}{http://www.research.att.com:80/cgi-bin/access.cgi/as/njas/sequences/eisA.cgi?Anum=A082161}
in \cite{oeis}.

\vspace*{10mm}

{\Large \textbf{4 \quad  Bijection from Paths to Automata}  }

In this section we exhibit a bijection from 
$\C_{k}^{\textrm{\raisebox{1ex}{$*$}}}(n,n)$ to 
$ \c_{k+1}(n) $. Using the illustrated path as a working example 
with $k=2$ and $n=4$,
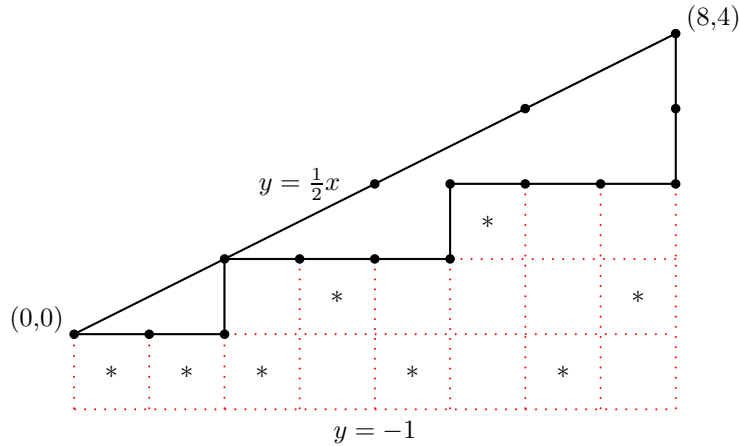
\begin{center} 
\begin{pspicture}(-4,0)(4,5.5)

\psline(-4,1)(-2,1)(-2,2)(1,2)(1,3)(4,3)(4,4)(4,5)
\psline(-4,1)(4,5)
\psline[linecolor=red,linestyle=dotted](-4,0)(4,0)
\psline[linecolor=red,linestyle=dotted](-2,1)(4,1)
\psline[linecolor=red,linestyle=dotted](1,2)(4,2)
\psline[linecolor=red,linestyle=dotted](-4,0)(-4,1)
\psline[linecolor=red,linestyle=dotted](-3,0)(-3,1)
\psline[linecolor=red,linestyle=dotted](-2,0)(-2,1)
\psline[linecolor=red,linestyle=dotted](-1,0)(-1,2)
\psline[linecolor=red,linestyle=dotted](0,0)(0,2)
\psline[linecolor=red,linestyle=dotted](1,0)(1,2)
\psline[linecolor=red,linestyle=dotted](2,0)(2,3)
\psline[linecolor=red,linestyle=dotted](3,0)(3,3)
\psline[linecolor=red,linestyle=dotted](4,0)(4,3)

\psdots(-4,1)(-3,1)(-2,1)(-2,2)(-1,2)(0,2)(1,2)(1,3)(2,3)(3,3)(4,3)(4,4)
\psdots(0,3)(2,4)(4,5)

\rput(-3.5,0.5){\textrm{{\footnotesize $*$}}}
\rput(-2.5,0.5){\textrm{{\footnotesize $*$}}}
\rput(-1.5,0.5){\textrm{{\footnotesize $*$}}}
\rput(-0.5,1.5){\textrm{{\footnotesize $*$}}}
\rput(0.5,0.5){\textrm{{\footnotesize $*$}}}
\rput(1.5,2.5){\textrm{{\footnotesize $*$}}}
\rput(2.5,0.5){\textrm{{\footnotesize $*$}}}
\rput(3.5,1.5){\textrm{{\footnotesize $*$}}}

\rput(-4.5,1.2){\textrm{{\footnotesize (0,0)}}}
\rput(4.5,5.2){\textrm{{\footnotesize (8,4)}}}
\rput(0,-0.3){\textrm{{\footnotesize $y=-1$}}}
\rput(-1,3.0){\textrm{{\footnotesize $y=\frac{1}{2}x$}}}

\end{pspicture}
\end{center} 
first construct the top row of a two-line representation consisting of 
$k+1$ each 1s,\,2s,\, \ldots,\,$n$\,s and number them left to right:
\[
\left(
\begin{array}{cccccccccccc}
    \overset{1}{1} &  \overset{2}{1}&  \overset{3}{1} &  \overset{4}{2} &  
    \overset{5}{2} &  \overset{6}{2} &  \overset{7}{3} &  
    \overset{8}{3}
    &  \overset{9}{3} &  \overset{10}{4} &  \overset{11}{4} &  
    \overset{12}{4}  \\
    ¥ & ¥ & ¥ & ¥ & ¥ & ¥ & ¥ & ¥ & ¥ & ¥ & ¥ & ¥
\end{array}\right).
\]
The last step in the path is necessarily an $N$ step. For the second 
last, third last,\ldots $N$ steps in the path, count the number of 
steps following it. This gives a sequence $i_{1},\,i_{2}, 
\ldots,i_{n-1}$ satisfying $1\le i_{1}<i_{2}<\ldots<i_{n-1}$ and 
$i_{j} \le (k+1)j$ for all $j$. Circle the positions $i_{1},\,i_{2}, 
\ldots,i_{n-1}$ in the two-line representation and then insert (in 
boldface) $2,3,\ldots,n$ in the second row in the circled positions:
\[
\left(
\begin{array}{cccccccccccc}
    \overset{\hspace*{-3mm}\textrm{\normalsize{ $\bigcirc$}} \hspace*{-3mm}1  }{1} &  \overset{2}{1}&  \overset{3}{1} &  \overset{4}{2} &  
    \overset{\hspace*{-3mm}\textrm{\normalsize{ $\bigcirc$}} 
    \hspace*{-3mm}5  }{2} &  \overset{6}{2} &  \overset{7}{3} &  
    \overset{8}{3}
    &  \overset{\hspace*{-3mm}\textrm{\normalsize{ $\bigcirc$}} 
    \hspace*{-3mm}9  }{3} &  \overset{10}{4} &  \overset{11}{4} &  
    \overset{12}{4}  \\
    \mathbf{2} & ¥ & ¥ & ¥ & \mathbf{3} & ¥ & ¥ & ¥ & \mathbf{4} & ¥ & ¥ & ¥
\end{array}\right).
\]
These will be the last occurrences of $2,3,\ldots,n$ in the second row. 
Working from the last column in the path back to the first, fill in 
the blanks in the second row left to right as follows. Count the number of 
squares from the $*$ up to the path (including the $*$ square) and 
add this number to the nearest boldface number to the left of the 
current blank entry (if there are  no boldface numbers to the left, 
add this number to 1) and insert the result in the current blank 
square. In the example the numbers of squares are 2,3,1,2,1,2,1,1 
yielding
\[
\left(
\begin{array}{cccccccccccc}
    \overset{\hspace*{-3mm}\textrm{\normalsize{ $\bigcirc$}} \hspace*{-3mm}1  }{1} &  \overset{2}{1}&  \overset{3}{1} &  \overset{4}{2} &  
    \overset{\hspace*{-3mm}\textrm{\normalsize{ $\bigcirc$}} 
    \hspace*{-3mm}5  }{2} &  \overset{6}{2} &  \overset{7}{3} &  
    \overset{8}{3}
    &  \overset{\hspace*{-3mm}\textrm{\normalsize{ $\bigcirc$}} 
    \hspace*{-3mm}9  }{3} &  \overset{10}{4} &  \overset{11}{4} &  
    \overset{12}{4}  \\
    \mathbf{2} & 4 & 5 & 3 & \mathbf{3} & 5 & 4 & 5 & \mathbf{4} & 5 & 5 & 
\end{array}\right).
\]
This will fill all blank entries except the last. Note that $*$\,s in 
the bottom row correspond to sink (that is, $n+1$) labels in the 
second row. Finally, insert $n+1$ into the last remaining blank space 
to give the image automaton: 
\[
\left(
\begin{array}{cccccccccccc}
   1 &  1&  1 & 2 &  
    2 & 2 &  3 &  3 &
    3 & 4 & 4 &  4  \\
    2 & 4 & 5 & 3 & 3 & 5 & 4 & 5 & 4 & 5 & 5 & 5
\end{array}\right).
\]
This process is fully reversible and the map is a bijection. \qed

\vspace*{10mm}

{\Large \textbf{5 \quad  Evaluation of det\,$\mathbf{A_{k}(n) }$  } }

For simplicity, we treat the case $k=1$, leaving the generalization to 
arbitrary $k$ as a not-too-difficult exercise for the interested reader. 
Write $A(n)$ for $A_{1}(n)$.
Thus $A(n)=\left( \StirlingCycle{\,i+1}{\,2i-j}\right)_{1\le i,j \le n}$.
From the definition of $\det 
A(n)$ as a sum of signed products, we  \linebreak \\[-1.1em]
show that $\det A(n)$ is 
the total weight of certain lists of permutations, each list 
carrying weight $\pm 1$. Then a weight-reversing involution cancels all $-1$ 
weights and reduces the problem to counting the surviving lists. These 
surviving lists are essentially the codes for paths in
$\C_{1}^{\textrm{\raisebox{1ex}{$*$}}}(n,p)$, and the Main Theorem 
follows from \S 4.

To describe the permutations giving a nonzero contribution to $\det 
A(n)=\sum_{\si}\textrm{sgn}\,\si\times \prod_{i=1}^{n}a_{i,\si(i)}$,
define the \emph{code} of a permutation $\si$ on $[n]$ to be the 
list $\textrm{\textbf{c}}=(c_{i})_{i=1}^{n}$ with 
$c_{i}=\si(i)-(i-1)$. Since the $(i,j)$ entry of $A(n),\ 
\StirlingCycle{i+1}{2i-j}$, is 0 unless $j\ge i-1$, we must have 
$\si(i)\ge i-1$ for all $i$. It is well known 
that there are $2^{n-1}$ such permutations, corresponding to 
compositions of $n$, with codes characterized by the following four 
conditions: (i) $c_{i}\ge 
0$ for all $i$, (ii) $c_{1}\ge 1$, (iii) each $c_{i}\ge 1$ is 
immediately followed by $c_{i}-1$ zeros in the list, (iv) 
$\sum_{i=1}^{n}c_{i}=n$. Let us call such a list a \emph{padded 
composition} of $n$: deleting the zeros is a bijection to ordinary 
compositions of $n$. For example, $(3,0,0,1,2,0)$ is a padded composition of 6. 
For a permutation $\si$ with padded composition code \textbf{c}, the 
nonzero entries in \textbf{c} give the cycle lengths of $\si$. Hence 
sgn\,$\si$, which is the parity of ``$n-\#\,$cycles in $\si$'', is 
given by $(-1)^{\textrm{\#\,0s in \textbf{c}}}$. 

We have $\det A(n) =\sum_{\si}\textrm{sgn}\,\si\,\prod_{i=1}^{n}a_{i,\si(i)} = 
\sum_{\si}\textrm{sgn}\,\si\,\prod_{i=1}^{n}\StirlingCycle{i+1}{2i-\si(i)} 
$, and so
\begin{equation}
    \det A(n)=\sum_{\textrm{\textbf{c}}} (-1)^{\textrm{\#\,0s in 
\textbf{c}}} \prod_{i=1}^{n}\StirlingCycle{i+1}{i+1-c_{i}}
    \label{eq:1}
\end{equation}
where the sum is restricted to padded compositions \textbf{c} of $n$ with 
$c_{i}\le i$ for all $i$
(\htmladdnormallink{A002083}{http://www.research.att.com:80/cgi-bin/access.cgi/as/njas/sequences/eisA.cgi?Anum=A002083})
because $\StirlingCycle{i+1}{i+1-c_{i}}=0$ unless $c_{i}\le i$.

Henceforth, let us write all permutations in standard cycle form 
whereby the smallest entry occurs first in each cycle and these 
smallest  entries increase left to right. Thus, with dashes 
separating cycles, 154-2-36 is the standard cycle form of the 
permutation $\left(
\begin{smallmatrix}
    1  & 2 & 3 & 4 & 5 & 6 \\
    5  & 2 & 6 & 1 & 4 & 3  
\end{smallmatrix}\right)$. We define a \emph{\ins} entry to be one that does not 
start a cycle. Thus the preceding permutation has 3 \ins entries: 
5,4,6. Note that the 
number of \ins entries is 0 only for the identity permutation. We 
denote an identity permutation (of any size) by $\ep$.

By definition of Stirling cycle number, the product in (\ref{eq:1}) counts lists 
$(\pi_{i})_{i=1}^{n}$ of permutations where $\pi_{i}$ is a 
permutation on $[i+1]$ with $i+1-c_{i}$ cycles, equivalently, with 
$c_{i}\le i$ \ins entries. 
So define $\el_{n}$ to be the set all lists of permutations 
$\pi=(\pi_{i})_{i=1}^{n}$ where $\pi_{i}$ is a 
permutation on $[i+1]$, \#\,\ins entries in $\pi_{i}$ is $\le i$, 
$\pi_{1}$ is the transposition (1,2), each nonidentity permutation $\pi_{i}$ 
is immediately followed by $c_{i}-1$ $\ep$'s where 
$c_{i}\ge 1$ is the number of \ins entries in $\pi_{i}$ (so the total 
number of \ins entries is $n$). Assign a weight to $\pi \in \el_{n}$ 
by wt$(\pi)=(-1)^{\textrm{\#\,$\ep$'s in $\pi$}}$.
Then 
\[
\det A(n) = \sum_{\pi \in \el_{n}}\textrm{wt}(\pi).
\]

We now define a weight-reversing involution on (most of) $\el_{n}$. 
Given $\mathbf{\pi} \in \el_{n}$, scan the list of its component 
permutations $\pi_{1}=(1,2),\pi_{2},\pi_{3},\ldots$ left to right. Stop 
at the first one that either (i) has more than one \ins entry, or (ii) 
has only one \ins entry, $b$ say, and $b>$ maximum \ins entry $m$ of the 
next permutation in the list. Say $\pi_{k}$ is the permutation where 
we stop.

In case (i) decrement (i.e. decrease by 1) the number of $\ep$'s in the list by 
splitting $\pi_{k}$ into two nonidentity permutations as follows. Let 
$m$ be the largest \ins entry of $\pi_{k}$ and let $\ell$ be 
its predecessor. Replace $\pi_{k}$ and its successor in the 
list (necessarily an $\ep$) by the following 
two permutations: first the transposition $(\ell,m)$ and second the 
permutation obtained from $\pi_{k}$ by erasing $m$ from its cycle and 
turning it into a singleton. Here are two examples of this case 
(recall permutations are in standard cycle form and, for clarity,
singleton cycles are not shown).
\[
\begin{array}{|c||c|c|c|c|c|c|}\hline 
    i & 1 & 2 & 3 & 4 & 5 & 6   \\ \hline
    \pi_{i} &   12 & 13 & 23 & 14\textrm{-}253 & \ep & \ep  \\ \hline
\end{array} \ \to \ 
\begin{array}{|c||c|c|c|c|c|c|}\hline
   i & 1 & 2 & 3 & 4 & 5 & 6   \\ \hline
    \pi_{i} &  12 & 13 & 23 & 25 & 14\textrm{-}23 & \ep   \\ \hline
\end{array}
\]
and
\[
\begin{array}{|c||c|c|c|c|c|c|}\hline 
    i & 1 & 2 & 3 & 4 & 5 & 6  \\ \hline
    \pi_{i} &   12 & 23 & 14 & 13\textrm{-}24 & \ep & 23 \\ \hline
\end{array} \ \to \ 
\begin{array}{|c||c|c|c|c|c|c|}\hline
   i & 1 & 2 & 3 & 4 & 5 & 6   \\ \hline
    \pi_{i} &  12 & 23 & 14 & 24 & 13 & 23 \\ \hline
\end{array}
\]

The reader may readily check that this sends case (i) to case (ii).

In case (ii), $\pi_{k}$ is a 
transposition $(a,b)$ with $b>$ maximum \ins entry $m$ of 
$\pi_{k+1}$. In this case, increment the number of $\ep$'s in the list by 
combining $\pi_{k}$ and $\pi_{k+1}$ into a single permutation 
followed by an $\ep$: in $\pi_{k+1}, b$ is a singleton; delete this singleton and 
insert $b$ immediately after $a$ in $\pi_{k+1}$ (in the same cycle). The 
reader may check that this reverses the result in the two examples above 
and, in general, sends case (ii) to case (i). Since the map alters the 
number of $\ep$'s in the list by 1, 
it is clearly weight-reversing. The map fails only for lists 
that both consist entirely of transpositions 
and have the form 
\[
(a_{1},b_{1}),\ (a_{2},b_{2}),\ \ldots,\ (a_{n},b_{n})\quad 
\textrm{with }b_{1}\le b_{2}\le \ldots \le b_{n}.
\]
Such lists have weight 1. Hence $\det A(n)$ is the number of lists 
$\big( (a_{i},b_{i})\big)_{i=1}^{n}$ satisfying $1 \le a_{i}<b_{i}\le 
i+1$ for $1 \le i \le n$, and $b_{1} \le b_{2} \le \ldots \le b_{n}$.
After subtracting 1 from each $b_{i}$, these lists code the paths 
in $\C_{1}^{\textrm{\raisebox{1ex}{$*$}}}(n,n)$ and, using \S 4, 
$\det A(n) = \v\, \C_{1}^{\textrm{\raisebox{1ex}{$*$}}}(n,n) \,\v =\v\, 
\c_{2}(n) \,\v$.

\end{document}